\documentclass[12pt]{article}

\usepackage{amsmath,amssymb,amsfonts,amsthm}
\usepackage{graphicx}
\usepackage{cite}
\usepackage[all]{xy}
\usepackage[toc,page]{appendix}
\usepackage{hyperref}
\hypersetup{colorlinks=true,  citecolor=red, linkcolor=blue}

\allowdisplaybreaks[3]

\newcounter{propositiona}
\newcommand{\propositiona}[1]{\refstepcounter{propositiona}
\noindent
\textbf{Proposition \thepropositiona.}\, {\it #1}}
\newcounter{definitiona}
\newcommand{\definitiona}[1]{\refstepcounter{definitiona}
\noindent
\textbf{Definition \thedefinitiona.}\, #1}
\newcounter{remarka}
\newcommand{\remarka}[1]{\refstepcounter{remarka}
\noindent
\textbf{Remark \theremarka.}\, #1}
\newcounter{examplea}
\newcommand{\examplea}[1]{\refstepcounter{examplea}
\noindent
\textbf{Example \theexamplea.}\, #1}
\newcounter{lemmaa}
\newcommand{\lemmaa}[1]{\refstepcounter{lemmaa}
\noindent
\textbf{Lemma \thelemmaa.}\, {\it #1}}
\newcounter{theorema}
\newcommand{\theorema}[1]{\refstepcounter{theorema}
\noindent
\textbf{Theorem\, \thetheorema.}\, {\it #1}}
\newcounter{corollarya}

\newcommand{\ord}{\mathrm{ord}\hspace{0.25ex} }

\title{A non-trivial conservation law with a trivial characteristic}

\author{ 
Kostya Druzhkov\hspace{0.1ex}\footnotemark[1] \vspace{0.5cm}\\
\small \emph{Department of Mathematics and Statistics, University of Saskatchewan, Saskatoon, Canada}\vspace{0.2cm}\\
}

\begin{document}

\footnotetext[1]{Electronic mail: konstantin.druzhkov@gmail.com}

\maketitle \numberwithin{equation}{section}

\vspace{-3ex}

\begin{abstract}
We show that the conservation law of the overdetermined system $u_t - 4u_x^3 - u_{xxx} = 0$, $u_y = 0$, associated with the characteristic $(u_{xy}, 0)$, is non-trivial despite the characteristic vanishing on the system.
\end{abstract}


\section{Introduction}

The notion of a conservation law requires a context—an analogue of phase space. One possible perspective is that such analogues are provided by infinitely prolonged differential equations, whose intrinsic geometry determines what should be considered a trivial conservation law. Historically, this perspective did not emerge immediately. At the time of the publication of Noether's seminal paper~\cite{Noether}, conservation laws were apparently understood in terms of divergence equations, without explicitly specifying what constituted trivial conservation laws. A proper understanding of the triviality of non-topological conservation laws arose later in~\cite{Steud}. A systematic approach based on the intrinsic geometry of differential equations emerged with the introduction of
the Vinogradov $\mathcal{C}$-spectral sequence~\cite{Vinold}.

From a practical point of view, the problem of the non-triviality of proper (non-topological) conservation laws turns out to be straightforward, for instance, for systems in Cauchy-Kovalevskaya form. It ultimately reduces to checking that a characteristic of a given conservation law is non-trivial, i.e., does not vanish on the corresponding system of equations~\cite{Alonso}. The following question arises: do there exist systems that admit proper conservation laws associated with trivial characteristics? This question was raised in~\cite{Olver0} for the Lagrangian case. However, it seems that since then, no such examples have been known, nor has their existence been confirmed, even for non-Lagrangian systems.

From a geometric point of view, it might be more meaningful to examine the existence of systems whose $\mathcal{C}$-spectral sequences have non-trivial, proper elements in their groups
$E_2^{\hspace{0.1ex} 0, \hspace{0.2ex} n-1}$ ($n$~denotes the number of independent variables). In general, elements of $E_2^{\hspace{0.1ex} 0, \hspace{0.2ex} n-1}$ are conservation laws that can be interpreted as rigid characteristic classes related to deformations of solutions with topologically non-trivial domains~\cite{Tsuj}. In this paper, we do not address topological questions but show that the group $E_2^{\hspace{0.1ex} 0, \hspace{0.2ex} 2}$ of the overdetermined system
\begin{align}
u_t - 4u_x^3 - u_{xxx} = 0\,,\qquad u_y = 0
\label{Intromain}
\end{align}
has a non-trivial element—the conservation law determined by
\begin{align*}
u_x^4\hspace{0.2ex} dt\wedge dx\,.
\end{align*}
One of its characteristics is $(u_{xy}, 0)$.
Even though the example is somewhat degenerate, its apparent simplicity is surprisingly misleading because, in fact, the group $E_2^{\hspace{0.1ex} 0, \hspace{0.2ex} 2}$ of~\eqref{Intromain} is isomorphic to the group $E_2^{\hspace{0.1ex} 2, \hspace{0.2ex} 1}$ of the potential mKdV equation
\begin{align*}
u_t - 4u_x^3 - u_{xxx} = 0\,,
\end{align*}
where $d_1$-closed elements of $E_1^{\hspace{0.1ex} 2, \hspace{0.2ex} 1}$ are its local presymplectic structures.

The paper is organized as follows. Section~\ref{Section2} introduces notation and provides basic facts from the geometry of differential equations. Section~\ref{Section3} is devoted to investigating a presymplectic structure of the potential mKdV equation. We demonstrate that the presymplectic operator given by the restriction of $D_x$ to this equation does not originate from its cosymmetries, and hence, leads to a non-trivial element of its group $E_2^{\hspace{0.1ex} 2, \hspace{0.2ex} 1}$. In Section~\ref{Section4}, we show that elements of the group $E_2^{\hspace{0.1ex} 2, \hspace{0.2ex} n-1}$ of a system of differential equations can be promoted to non-trivial conservation laws of some closely related system. Such conservation laws are associated with trivial characteristics. We provide details on the example announced above. Although the potential mKdV is integrable, Theorem~\ref{TheorLagrExa} shows that its integrability is not essential in the context of the questions under consideration.

We use Einstein's summation convention throughout this paper and consider only smooth
functions of the class $C^{\infty}$.

\section{\label{Section2}Basic notation and concepts}

Let us introduce notation and recall basic facts from the geometry of differential equations. More details can be found in~\cite{VinKr, KraVer1}.

\subsection{Jets}

Let $\pi\colon E^{n+m}\to M^n$ be a locally trivial smooth vector
bundle over a smooth manifold $M^n$. The bundle $\pi$ gives rise to the corresponding jet bundles $\pi_k\colon J^k(\pi)\to M$,
\begin{align*}
\xymatrix{
\ldots \ar[r] & J^3(\pi) \ar[r]^-{\pi_{3, 2}} & J^2(\pi) \ar[r]^-{\pi_{2, 1}} & J^1(\pi) \ar[r]^-{\pi_{1, 0}} & J^0(\pi) = E \ar[r]^-\pi & M
}
\end{align*}
and the inverse limit $J^{\infty}(\pi)$ arising with the natural projections 
$\pi_{\infty}\colon J^{\infty}(\pi)\to M$ and $\pi_{\infty,\hspace{0.2ex} k}\colon J^{\infty}(\pi)\to J^k(\pi)$. 
Denote by $\mathcal{F}(\pi)$ the algebra of smooth functions on~$J^{\infty}(\pi)$,
\begin{align*}
\mathcal{F}(\pi) = \bigcup_{k\geqslant 0} \pi_{\infty,\hspace{0.2ex} k}^{\hspace{0.1ex} *}\, C^{\infty}(J^k(\pi))\,.
\end{align*}

Smooth sections of the pullback bundles $\pi_k^*(\pi)\colon \pi_k^*(E)\to J^k(\pi)$ can be interpreted as sections of $\pi_{\infty}^*(\pi)\colon \pi_{\infty}^*(E) \to J^{\infty}(\pi)$. Denote by $\varkappa(\pi)$ the $\mathcal{F}(\pi)$-module $\Gamma(\pi_{\infty}^*(\pi))$ of such sections of $\pi_{\infty}^*(\pi)$.


\vspace{0.5ex}

\noindent
\textbf{Local coordinates.} Suppose $U\subset M$ is a coordinate neighborhood such that the bundle $\pi$
becomes trivial over $U$. Choose local coordinates $x^1$, \ldots, $x^n$ in $U$ and $u^1$, \ldots, $u^m$  
along the fibers of $\pi$ over $U$. Then a section $\sigma\in \Gamma(\pi)$ takes the form of a smooth vector function $u^i = \sigma^i(x^1, \ldots, x^n)$.
We denote by $u^i_{\alpha}$ the corresponding adapted local coordinates on $J^{\infty}(\pi)$, where $\alpha$ denotes a multi-index.
It is convenient to treat $\alpha$ as a formal sum of the form $\alpha = \alpha_1 x^1 + \ldots + \alpha_n x^n = \alpha_i x^i$, where all $\alpha_i$ are non-negative integers; $|\alpha| = \alpha_1 + \ldots + \alpha_n$. For the infinite jet $[\sigma]_{x_0}^{\infty}$ of a section $\sigma\in \Gamma(\pi)$ at a point $x_0\in U$, we have
\begin{align*}
x^i([\sigma]^{\infty}_{x_0}) = x^i_0\,,\qquad u^i_{\alpha}([\sigma]^{\infty}_{x_0}) = \dfrac{\partial^{|\alpha|} \sigma^i}{(\partial x^1)^{\alpha_1}\ldots (\partial x^n)^{\alpha_n}}(x_0)\,.
\end{align*}
In what follows, we consider only adapted local coordinates on $J^{\infty}(\pi)$. 

\vspace{0.5ex}

\noindent
\textbf{Cartan distribution.} The main structure on jet manifolds is the Cartan distribution.
Using adapted local coordinates on $J^{\infty}(\pi)$, one can introduce the total derivatives
$$
D_{x^i} = \partial_{x^i} + u^k_{\alpha + x^i}\hspace{0.15ex} \partial_{u^k_{\alpha}}\qquad\quad i = 1, \ldots, n.
$$
The planes of the Cartan distribution $\mathcal{C}$ on $J^{\infty}(\pi)$ are spanned by the total derivatives.

\vspace{0.5ex}
\noindent
\textbf{Cartan forms.} The Cartan distribution $\mathcal{C}$ determines the ideal $\mathcal{C}\Lambda^*(\pi)$
of the algebra
$$
\Lambda^*(\pi) = \bigcup_{k\geqslant 0} \pi_{\infty,\hspace{0.2ex} k}^{\hspace{0.1ex} *}\, \Lambda^*(J^k(\pi))
$$
of differential forms on $J^{\infty}(\pi)$.
The ideal $\mathcal{C}\Lambda^*(\pi)$ is generated by Cartan (or contact) forms, i.e., differential forms that vanish on any plane of the Cartan distribution $\mathcal{C}$.
A Cartan $1$-form $\omega\in\mathcal{C}\Lambda^1(\pi)$ can be written as a finite sum
$$
\omega = \omega_i^{\alpha}\theta^i_{\alpha}\,,\qquad\ \theta^i_{\alpha} = du^i_{\alpha} - u^i_{\alpha + x^k}dx^k
$$
in adapted local coordinates. The coefficients $\omega_i^{\alpha}$ are smooth functions of adapted coordinates.

\vspace{0.5ex}
\noindent
\textbf{Infinitesimal symmetries.} Each section $\varphi\in \varkappa(\pi)$ gives rise to the corresponding evolutionary vector field $E_{\varphi}$ on $J^{\infty}(\pi)$,
$$
E_{\varphi} = D_{\alpha}(\varphi^i)\partial_{u^i_{\alpha}}\,.
$$
Here $\varphi^1$, \ldots, $\varphi^m$ are components of $\varphi$, $D_{\alpha}$ denotes the composition $D_{x^1}^{\ \alpha_1}\circ\ldots\circ D_{x^n}^{\ \alpha_n}$.
Evolutionary vector fields are infinitesimal symmetries of $J^{\infty}(\pi)$. In particular, $\mathcal{L}_{E_{\varphi}}\,\mathcal{C}\Lambda^*(\pi)\subset \mathcal{C}\Lambda^*(\pi)$. Here $\mathcal{L}_{E_{\varphi}}$ is the corresponding Lie derivative. Elements of $\varkappa(\pi)$ are characteristics of symmetries of $J^{\infty}(\pi)$.

\vspace{0.5ex}
\noindent
\textbf{Horizontal forms.}
Cartan forms allow one to consider the modules of horizontal $k$-forms
$$
\Lambda^k_h(\pi) = \Lambda^k(\pi)/\mathcal{C}\Lambda^k(\pi)\,.
$$
The de Rham differential $d$ induces the horizontal differential $d_h\colon \Lambda^k_h(\pi)\to \Lambda^{k+1}_h(\pi)$.
The infinite jet bundle $\pi_{\infty}\colon J^{\infty}(\pi) \to M$ admits the decomposition
$$
\Lambda^1(\pi) = \mathcal{C}\Lambda^1(\pi) \oplus \mathcal{F}(\pi)\!\cdot\!\pi^*_{\infty}(\Lambda^1(M))\,.
$$
We identify the module of horizontal $k$-forms $\Lambda^k_h(\pi)$ with $\mathcal{F}(\pi)\cdot \pi^*_{\infty}(\Lambda^k(M))$.
In adapted local coordinates, elements of $\mathcal{F}(\pi)\cdot \pi^*_{\infty}(\Lambda^k(M))$ are generated by the differentials $dx^1, \ldots, dx^n$, while $d_h = dx^i\wedge D_{x^i}$. For example,
$$
d_h (\xi_j dx^j) = dx^i\wedge D_{x^i}(\xi_j) dx^j = D_{x^i}(\xi_j)dx^i\wedge dx^j\,.
$$

\noindent
\textbf{Euler operator.} Let $\widehat{\varkappa}(\pi)$ be the adjoint module
\begin{align*}
\widehat{\varkappa}(\pi) = \mathrm{Hom}_{\mathcal{F}(\pi)}(\varkappa(\pi), \Lambda^n_h(\pi))\,.
\end{align*}
Denote by $\mathrm{E}$ the Euler operator (variational derivative), $\mathrm{E}\colon \Lambda^n_h(\pi)\to \widehat{\varkappa}(\pi)$.
In adapted local coordinates, for $L = \lambda \, dx^1\wedge\ldots\wedge dx^n$ and $\varphi\in\varkappa(\pi)$, we have
\begin{align*}
&\mathrm{E}(L)\colon\varphi \mapsto \langle \mathrm{E}(L), \varphi\rangle = \dfrac{\delta \lambda}{\delta u^i}\,\varphi^i dx^1\wedge\ldots\wedge dx^n\,, \qquad \dfrac{\delta \lambda}{\delta u^i} = \sum_{\alpha} (-1)^{|\alpha|}D_{\alpha}\Big(\dfrac{\partial \lambda}{\partial u^i_{\alpha}}\Big)\,.
\end{align*}
Here $\langle \cdot, \cdot \rangle$ denotes the natural pairing between a module and its adjoint.

\subsection{Differential equations}

Let $\zeta \colon E_1\to M$ be a locally trivial smooth vector bundle over the same base as $\pi$. Smooth sections of the pullbacks $\pi^*_{r}(\zeta)$ determine a module of sections of the pullback $\pi^*_{\infty}(\zeta)\colon \pi^*_{\infty}(E_1)\to J^{\infty}(\pi)$. We denote it by $P(\pi)$. Any $F\in P(\pi)$ can be considered a (generally, nonlinear) differential operator $\Gamma(\pi)\to \Gamma(\zeta)$. Then $F = 0$ is a differential equation.
By its \emph{infinite prolongation} we mean the set of formal solutions $\mathcal{E}\subset J^{\infty}(\pi)$ defined by the infinite system of equations
\begin{align*}
\mathcal{E}\colon\qquad D_{\alpha}(F^i) = 0\,,\qquad |\alpha| \geqslant 0\,.
\end{align*}
Here $F^i$ are components of $F$ in adapted coordinates. We assume that $\pi_{\infty}(\mathcal{E}) = M$.

\vspace{1ex}
\remarka{We do not require that the number of equations of the form $F^i = 0$ coincide with the number of dependent variables $m$.}

\vspace{1ex}
\noindent
\textbf{Regularity assumptions.} We consider only systems that satisfy the following conditions.
\begin{enumerate}
  \item For any $\rho\in \{F = 0\}\subset J^{r}(\pi)$, the differentials $dF^i_\rho$ of the components $F^i$ are independent.
  \item A function $f\in \mathcal{F}(\pi)$ vanishes on $\mathcal{E}$ if and only if there exists a differential operator $\Delta\colon P(\pi)\to \mathcal{F}(\pi)$ of the form $\Delta_i^{\alpha} D_{\alpha}$ ($\mathcal{C}$-differential operator or total differential operator) such that $f = \Delta(F)$. Here, for some integer $k$, $|\alpha| \leqslant k$ for all components $\Delta_i^{\alpha}$. These components may depend on independent variables $x^i$, dependent variables $u^i$, and derivatives up to the $k$-th order.

\end{enumerate}

For simplicity, we assume that the de Rham cohomology groups $H^i_{dR}(\mathcal{E})$ are trivial for $i > 0$.

\vspace{0.5ex}
\noindent
\textbf{Functions.}
By $\mathcal{F}(\mathcal{E})$ we denote the algebra of smooth functions on $\mathcal{E}$,
$$
\mathcal{F}(\mathcal{E}) = \mathcal{F}(\pi)|_{\mathcal{E}} = \mathcal{F}(\pi)/I\,.
$$
Here $I$ denotes the ideal of the system $\mathcal{E}\subset J^{\infty}(\pi)$, $I = \{f\in \mathcal{F}(\pi)\, \colon\, f|_{\mathcal{E}} = 0\}$.

\vspace{0.5ex}
\noindent
\textbf{Cartan forms.} The module $\Lambda^1(\mathcal{E}) = \Lambda^1(\pi)|_{\mathcal{E}} = \Lambda^1(\pi)/(I\cdot \Lambda^1(\pi) + \mathcal{F}(\pi)\cdot dI)$ of differential $1$-forms on $\mathcal{E}$ produces the exterior algebra $\Lambda^*(\mathcal{E})$.
The Cartan distribution of $J^{\infty}(\pi)$ can be restricted to $\mathcal{E}$.
Similarly, the ideal $\mathcal{C}\Lambda^*(\mathcal{E})\subset \Lambda^*(\mathcal{E})$ is generated by differential forms that vanish on any plane of the Cartan distribution of $\mathcal{E}$.
Note that $\mathcal{C}\Lambda^*(\mathcal{E}) = \mathcal{C}\Lambda^*(\pi)|_{\mathcal{E}}$ due to the decomposition of $\Lambda^1(\pi)$.

\vspace{0.5ex}
\noindent
\textbf{Infinitesimal symmetries.} A \emph{symmetry} (more precisely, an infinitesimal symmetry) of an infinitely prolonged system of equations $\mathcal{E}$ is a vector field $X$ on $\mathcal{E}$ (a derivation of $\mathcal{F}(\mathcal{E})$) that preserves the Cartan distribution, i.e., $\mathcal{L}_X\, \mathcal{C}\Lambda^*(\mathcal{E})\subset \mathcal{C}\Lambda^*(\mathcal{E})$ for $\mathcal{L}_X = X \lrcorner \circ d + d \circ X \lrcorner\,$. Two symmetries are equivalent if they differ by a trivial symmetry, i.e., a vector field on $\mathcal{E}$ such that at each point of $\mathcal{E}$, its vector lies in the respective plane of the Cartan distribution. One can say that, locally, trivial symmetries are combinations of the total derivatives $\,\overline{\!D}_{x^i} = D_{x^i}|_{\mathcal{E}}$, $i = 1, \ldots, n$.

If $\varphi\in \varkappa(\pi)$ is a characteristic such that $E_{\varphi}$ is tangent to $\mathcal{E}$ (i.e., $E_{\varphi}(F)|_{\mathcal{E}} = 0$, or equivalently, $E_{\varphi}(I)\subset I$), then the restriction $E_{\varphi}|_{\mathcal{E}}\colon \mathcal{F}(\mathcal{E})\to \mathcal{F}(\mathcal{E})$ is a symmetry of $\mathcal{E}$. In other words, elements of the kernel of the linearization operator $l_{\mathcal{E}} = l_F|_{\mathcal{E}}\colon \varkappa(\mathcal{E})\to P(\mathcal{E})$ correspond to symmetries of $\mathcal{E}$. Here $l_F\colon \varkappa(\pi)\to P(\pi)$, $\varphi \mapsto E_{\varphi}(F)$, $l_F(\varphi)^i = E_{\varphi}(F^i)$, and
$$
\varkappa(\mathcal{E}) = \varkappa(\pi)|_{\mathcal{E}} = \varkappa(\pi)/I\cdot \varkappa(\pi)\,,\qquad P(\mathcal{E}) = P(\pi)|_{\mathcal{E}} = P(\pi)/I\cdot P(\pi)\,.
$$

If $\pi_{\infty,\, 0}(\mathcal{E}) = J^0(\pi)$, then for each symmetry $X$ of $\mathcal{E}\subset J^{\infty}(\pi)$, there exists $\varphi \in \varkappa(\pi)$ such that $X$ is equivalent to the restriction $E_{\varphi}|_{\mathcal{E}}$. In this case, symmetries of $\mathcal{E}\subset J^{\infty}(\pi)$ are in one-to-one correspondence with elements of $\ker l_{\mathcal{E}}$.

\vspace{0.5ex}
\noindent
\textbf{$\mathcal{C}$-spectral sequence.} For $p \geqslant 1$, the $p$-th power $\mathcal{C}^p\Lambda^*(\mathcal{E})$ of the ideal $\mathcal{C}\Lambda^*(\mathcal{E})$ is stable with respect to the de Rham differential $d$, i.e., $d(\mathcal{C}^p\Lambda^*(\mathcal{E})) \subset \mathcal{C}^p\Lambda^*(\mathcal{E})$.
Then the de Rham complex $\Lambda^{\bullet}(\mathcal{E})$ admits the filtration
$$
\Lambda^{\bullet}(\mathcal{E})\supset \mathcal{C}\Lambda^{\bullet}(\mathcal{E})\supset \mathcal{C}^2\Lambda^{\bullet}(\mathcal{E})\supset \mathcal{C}^3\Lambda^{\bullet}(\mathcal{E})\supset \ldots
$$
The corresponding spectral sequence $(E^{\hspace{0.1ex} p,\hspace{0.2ex} q}_r(\mathcal{E}), d^{\hspace{0.1ex} p,\hspace{0.2ex} q}_r)$ is the Vinogradov $\mathcal{C}$-spectral sequence~\cite{Vin, VinKr}.
Here $\mathcal{C}^{k+1}\Lambda^k(\mathcal{E}) = 0$, $\mathcal{C}^0\Lambda^q(\mathcal{E}) = \Lambda^q(\mathcal{E})$, $E^{\hspace{0.1ex} p, \hspace{0.2ex} q}_0(\mathcal{E}) = \mathcal{C}^p\Lambda^{p+q}(\mathcal{E})/\mathcal{C}^{p+1}\Lambda^{p+q}(\mathcal{E})$. All the differentials $d_r^{\hspace{0.1ex} p,\hspace{0.2ex} q}$ are induced by the de Rham differential $d$,
$$
d_r^{\hspace{0.1ex} p,\hspace{0.2ex} q}\colon E^{\hspace{0.1ex} p, \hspace{0.2ex} q}_r(\mathcal{E}) \to E^{\hspace{0.1ex} p+r, \hspace{0.2ex} q+1-r}_r(\mathcal{E})\,,\qquad E^{\hspace{0.1ex} p, \hspace{0.2ex} q}_{r+1}(\mathcal{E}) = \ker d_r^{\hspace{0.1ex} p,\hspace{0.2ex} q}/ \mathrm{im}\, d_r^{\hspace{0.1ex} p-r,\hspace{0.2ex} q+r-1}\qquad \text{for}\quad r\geqslant 0\,.
$$
We also use the notation $d_r$ where it does not lead to confusion. Besides, we put $\Lambda^k_h(\mathcal{E}) = E_0^{\hspace{0.1ex} 0, \hspace{0.2ex} k}(\mathcal{E})$. Then the horizontal differential $d_h\colon \Lambda^k_h(\mathcal{E})\to \Lambda^{k+1}_h(\mathcal{E})$ coincides with $d^{\hspace{0.2ex} 0, \hspace{0.2ex} k}_0$. The $\mathcal{C}$-spectral sequence of a differential equation converges to its de Rham cohomology.

$\mathcal{C}$-spectral sequence allows one to define variational $k$-forms of differential equations. A \emph{variational $k$-form} of $\mathcal{E}$ is an element of the group $E^{\hspace{0.1ex} k,\hspace{0.2ex} n-1}_1(\mathcal{E})$. The group $E^{\hspace{0.1ex} k,\hspace{0.2ex} n-1}_1(\mathcal{E})$ is canonically isomorphic to the group
\begin{align*}
\dfrac{\{\omega\in \mathcal{C}^k\Lambda^{k+n-1}(\mathcal{E})\hspace{0.1ex} \colon\ d\omega \in \mathcal{C}^{k+1}\Lambda^{k+n}(\mathcal{E})\}}{\mathcal{C}^{k+1}\Lambda^{k+n-1}(\mathcal{E}) + d(\mathcal{C}^k\Lambda^{k+n-2}(\mathcal{E}))}\,.
\end{align*}
Then any $\omega\in \mathcal{C}^k\Lambda^{k+n-1}(\mathcal{E})$ such that $d\omega \in \mathcal{C}^{k+1}\Lambda^{k+n}(\mathcal{E})$ gives rise to a variational $k$-form. For any symmetry $X$, the action of the Lie derivative $\mathcal{L}_X$ on elements of $E^{\hspace{0.1ex} p,\hspace{0.2ex} q}_r(\mathcal{E})$ is well-defined.

\vspace{0.5ex}

\noindent
\textbf{Conservation laws.} A \emph{conservation law} of $\mathcal{E}$ is a variational $0$-form, i.e., an element of the group $E^{\hspace{0.1ex} 0,\hspace{0.2ex} n-1}_1(\mathcal{E})$. Let us denote by $\widehat{P}(\pi)$ the adjoint module
\begin{align*}
\widehat{P}(\pi) = \mathrm{Hom}_{\mathcal{F}(\pi)}(P(\pi), \Lambda^n_h(\pi))\,.
\end{align*}
A \textit{characteristic} of a conservation law~\cite{Olver} is any homomorphism $Q\in \widehat{P}(\pi)$ such that
\begin{align*}
\langle Q, F \rangle = d_h \mu\,,\qquad \mu\in\Lambda_h^{n-1}(\pi)\,,
\end{align*}
where $\mu|_{\mathcal{E}}\in\Lambda_h^{n-1}(\mathcal{E}) = E_0^{\hspace{0.1ex} 0,\hspace{0.2ex} n-1}(\mathcal{E})$ represents the conservation law. The total divergence form of the conservation law is $D_{x^i}(\mu^i) = 0$, where $\mu = \mu^i \partial_{x^i} \lrcorner\, (dx^1\wedge \ldots\wedge dx^n)$.

\vspace{0.5ex}

\noindent
\textbf{Cosymmetries.} Elements of the kernel of the adjoint linearization operator $l_{\mathcal{E}}^{\, *}\colon \widehat{P}(\mathcal{E})\to \widehat{\varkappa}(\mathcal{E})$ are called \emph{cosymmetries} of $\mathcal{E}$. Here $l_{\mathcal{E}}^{\, *} = l_{F}^{\, *}|_{\mathcal{E}}$, $\widehat{P}(\mathcal{E}) = \widehat{P}(\pi)|_{\mathcal{E}}$, $\widehat{\varkappa}(\mathcal{E}) = \widehat{\varkappa}(\pi)|_{\mathcal{E}}$. The restriction $Q|_{\mathcal{E}}$ of a characteristic $Q$ of a conservation law is its cosymmetry. However, a more precise statement is that cosymmetries encode variational $1$-forms. Let $\psi\in \ker l_{\mathcal{E}}^{\, *}$ be a cosymmetry, and let $G\in \widehat{P}(\pi)$ be a homomorphism such that $G|_{\mathcal{E}} = \psi$. According to the Green formula (the definition of $l_F^{\,*}$), there exists a differential form $\chi\in \mathcal{C}\Lambda^1(\mathcal{\pi})\wedge \pi^*_{\infty}(\Lambda^{n-1}(M))$ such that
\begin{align}
\langle G, l_F(\varphi) \rangle = \langle l_F^{\,*}(G), \varphi \rangle + d_h (E_{\varphi} \lrcorner\, \chi) \qquad \text{for any}\quad \varphi\in \varkappa(\pi)\,.
\label{Cosymtoonef}
\end{align}
The differential form $\chi|_{\mathcal{E}}$ gives rise to a variational $1$-form of $\mathcal{E}$. This map from cosymmetries to variational $1$-forms is surjective for any system $\mathcal{E}$ (see Appendix~\ref{App:A}). The correspondence is an isomorphism if $\mathcal{E}$ is $\ell$-normal~\cite{VinKr} (i.e., such that the equation $\nabla \circ l_{\mathcal{E}} = 0$ for a total differential operator $\nabla\colon P(\mathcal{E})\to \mathcal{F}(\mathcal{E})$ has no nonzero solutions). The converse is also true (see Appendix~\ref{App:A}).
In the general case, a cosymmetry of a conservation law $\xi \in E^{\hspace{0.1ex} 0,\hspace{0.2ex} n-1}_1(\mathcal{E})$ is any cosymmetry that corresponds to $d_1\xi$.
A cosymmetry is \textit{trivial} if the corresponding variational $1$-form is zero.


\vspace{0.5ex}
\noindent
\textbf{Presymplectic structures.} A \emph{presymplectic structure} of $\mathcal{E}$ is a $d_1$-closed variational $2$-form, i.e., an element of the kernel of the differential
$$
d_1^{\hspace{0.2ex} 2,\hspace{0.2ex} n-1}\colon E^{\,2,\,n-1}_1(\mathcal{E})\to E^{\,3,\,n-1}_1(\mathcal{E}).
$$

\vspace{0.5ex}

\noindent
\textbf{Internal Lagrangians.} Elements of the group 
\begin{align*}
\widetilde{E}_1^{\hspace{0.1ex} 0,\hspace{0.2ex} n-1}(\mathcal{E}) = \dfrac{\{l\in \Lambda^{n}(\mathcal{E})\hspace{0.1ex} \colon \ dl\in \mathcal{C}^2\Lambda^{n+1}(\mathcal{E})\}}{\mathcal{C}^2\Lambda^{n}(\mathcal{E}) + d(\Lambda^{n-1}(\mathcal{E}))}
\end{align*}
are called \emph{internal Lagrangians}~\cite{Druzhkov2, Druzhkov3} of $\mathcal{E}$. The de Rham differential $d$ induces the differential 
$$
\tilde{d}^{\hspace{0.3ex} 0,\hspace{0.2ex} n-1}_1\colon \widetilde{E}_1^{\hspace{0.1ex} 0,\hspace{0.2ex} n-1}(\mathcal{E})\to E^{\hspace{0.1ex}2,\hspace{0.2ex} n-1}_1(\mathcal{E})
$$
mapping them to presymplectic structures of $\mathcal{E}$.

If the variational derivative $\mathrm{E}(L)$ of a horizontal $n$-form $L\in \Lambda^n_h(\pi)$ vanishes on $\mathcal{E}$, the horizontal cohomology class $L + d_h \Lambda^{n-1}_h(\pi)$ produces a unique internal Lagrangian of $\mathcal{E}$ via the Noether formula:
there exists a differential form $\omega_L\in \mathcal{C}\Lambda^1(\mathcal{\pi})\wedge \pi^*_{\infty}(\Lambda^{n-1}(M))$ such that for any $\varphi\in \varkappa(\pi)$,
\begin{align}
\mathcal{L}_{E_\varphi}(L) = \langle \mathrm{E}(L), \varphi \rangle + d_h (E_\varphi \lrcorner\, \omega_L)\,.
\label{Noethiden}
\end{align}
Then the differential $n$-form $l = (L + \omega_L)|_{\mathcal{E}}$ represents the corresponding internal Lagrangian, while $dl$ gives rise to the presymplectic structure. The formula~\eqref{Noethiden} expresses integration by parts.


\subsection{Some computational formulae}

From a computational point of view, below, we deal only with a scalar $1+1$ evolution equation. Evolution systems are convenient because their presymplectic structures can be identified with operators of a certain form~\cite{Druzhkov1} rather than equivalence classes of operators~\cite{VinKr}. If
$$
F^i = u^i_t - \Phi^i(t, x^2, \ldots, x^n, u^1, \ldots, u^n, u^1_{x^2}, \ldots)\,,\qquad x^1 = t\,,
$$
then all the variables except for the ones of the form $u^i_{t + \alpha}$ can be treated as coordinates on $\mathcal{E}$. In this case, $\mathcal{F}(\mathcal{E})\subset \mathcal{F}(\pi)$, $\varkappa(\mathcal{E})\subset \varkappa(\pi)$, $\widehat{P}(\mathcal{E})\subset \widehat{P}(\pi)$, and one has the following simplifications.
\begin{enumerate}
  \item For $\varphi\in \ker l_{\mathcal{E}}\subset \varkappa(\pi)$ and $\psi\in \ker l_{\mathcal{E}}^{\,*}\subset \widehat{P}(\pi)$,
\begin{align*}
l_F(\varphi) = l_{\varphi}(F)\,,\qquad l_F^{\,*}(\psi) = -l_{\psi}(F)\,.
\end{align*}
  \item All presymplectic structures are extendable~\cite{Druzhkov2}, i.e., originate from Lagrangians. If the variational derivative $\mathrm{E}(L)$ of a Lagrangian $L\in \Lambda_h^n(\pi)$ vanishes on $\mathcal{E}$, then $\mathrm{E}(L) = A(F)$ for some total differential operator $A\colon P(\pi)\to \widehat{\varkappa}(\pi)$. The presymplectic structure produced by $L$ corresponds to a unique operator of the form $A^*|_{\mathcal{E}} + \square\circ l_{\mathcal{E}}$ that does not involve $D_t|_{\mathcal{E}}$.
  \item If a variational $1$-form $\omega_{\psi}$ corresponds to a cosymmetry $\psi$, then the presymplectic structure $d_1\omega_{\psi}$ corresponds to the presymplectic operator $l_{\psi} - l_{\psi}^{\,*}$\hspace{0.15ex}, which does not involve $D_t|_{\mathcal{E}}$.
  \item If $X$ is a symmetry equivalent to $E_{\varphi}|_{\mathcal{E}}$, then, according to the relation between cosymmetries and variational $1$-forms based on~\eqref{Cosymtoonef}, the cosymmetry of the conservation law $X\lrcorner\, \omega_{\psi}$ is 
$$
l_{\varphi}^{\,*}(\psi) + l_{\psi}^{\,*}(\varphi)\,.
$$
  \item The cosymmetry of $\mathcal{L}_X \omega_{\psi} = X\lrcorner\, d_1\omega_{\psi} + d_1 (X\lrcorner\, \omega_{\psi})$ is $l_{\psi}(\varphi) - l_{\psi}^{\,*}(\varphi) + l_{\varphi}^{\,*}(\psi) + l_{\psi}^{\,*}(\varphi)$, i.e.,
$$
E_{\varphi}(\psi) + l_{\varphi}^{\,*}(\psi)\,.
$$

\end{enumerate}

\section{\label{Section3}A presymplectic structure that is not $d_1$-exact}

Any $\ell$-normal system of differential equations has trivial group $E_2^{\hspace{0.1ex} 1, \hspace{0.2ex} n-1}$ of the Vinogradov $\mathcal{C}$-spectral sequence~\cite{VinKr}. Then cosymmetries of $\ell$-normal systems produce either non-trivial conservation laws or non-trivial presymplectic structures.
More generally, cosymmetries of an arbitrary system of differential equations produce either its conservation laws or non-trivial internal Lagrangians (see~\cite{Druzhkov3}, Section~4). For some systems, all their local presymplectic structures originate from cosymmetries. This is the case, for example, for any system that admits a symmetry group scaling each dependent variable by a factor with a positive exponent (not necessarily the same for all dependent variables) while preserving all independent variables (such systems are called conic equations in~\cite{Vin}; see Proposition on p.~111). In particular, all linear systems are of this type. We provide an example of a presymplectic structure that cannot be obtained through cosymmetries.

Let us consider the infinite prolongation $\mathcal{E}$ of the potential mKdV equation
\begin{align}
\label{almostpKdV}
u_t = 4u_x^3 + u_{xxx}\,.
\end{align}
Here $n=2$, $m=1$, $\pi\colon \mathbb{R}\times \mathbb{R}^2\to \mathbb{R}^2$ is the projection onto the second factor.

The stationary action principle for the Lagrangian
\begin{align}
L = \lambda\, dt\wedge dx\,,\qquad \lambda = \dfrac{u_t u_x}{2} - u_x^4 + \dfrac{u_{xx}^2}{2}
\label{almostpKdVLagr}
\end{align}
leads to a differential consequence of~\eqref{almostpKdV}. Then $L$ gives rise to a unique internal Lagrangian of $\mathcal{E}$ and to the corresponding presymplectic structure $\Omega\in E_1^{\hspace{0.1ex} 2,\hspace{0.2ex} 1}(\mathcal{E})$, $d_1\Omega = 0$.

As coordinates on $\mathcal{E}$, we take $t$, $x$, $u_0 = u$, $u_1 = u_x$, $u_2 = u_{xx}$, $u_3 = u_{xxx}$, $u_4 = u_{xxxx}$, $\ldots$ Each function from the algebra $\mathcal{F}(\mathcal{E})$ depends on a finite number of the coordinates. It is convenient to introduce the concept of order of a function, which can be a non-negative integer or $-\infty$. It is simply the highest order of the derivatives $u_i$ among its arguments.

\vspace{1ex}

\definitiona{A function $f\in\mathcal{F}(\mathcal{E})$ has \textit{order $k$} if $\partial_{u_k}\hspace{0.1ex} f$ is not identically zero, whereas $\partial_{u_{k+h}} f \equiv 0$ for any integer $h\geqslant 1$. In this case, we use the notation $\ord f = k$. If $f$ does not depend on the coordinates of the form $u_i$ ($i = 0, 1, \ldots$), then $\ord f = -\infty$.
}

\vspace{1ex}
For~\eqref{almostpKdV}, one has $\widehat{P}(\mathcal{E}) = \Lambda^2_h(\mathcal{E}) = \widehat{\varkappa}(\mathcal{E}) \simeq \varkappa(\mathcal{E}) = P(\mathcal{E}) = \mathcal{F}(\mathcal{E})$. We identify $\widehat{P}(\mathcal{E})$ with~$\mathcal{F}(\mathcal{E})$ taking the $dt\wedge dx$-components of elements of $\widehat{P}(\mathcal{E}) = \Lambda^2_h(\mathcal{E})$ (by applying $\,\overline{\!D}_x\lrcorner \ \overline{\!D}_t\lrcorner\, $ to them). Here the total derivatives span the Cartan distribution, 
\begin{align*}
\,\overline{\!D}_x = \partial_x + u_{i+1}\partial_{u_{i}}\,,\qquad 
\,\overline{\!D}_t = \partial_t + \,\overline{\!D}_x^{\,i}(4u_1^3 + u_{3})\partial_{u_{i}}\,.
\end{align*}
Let us note that the operator $\,\overline{\!D}_x$ increases the order of any function by $1$, and its kernel consists of functions of $t$.

The presymplectic structure $\Omega$ corresponds to the presymplectic operator
\begin{align*}
\Delta = \,\overline{\!D}_x\,.
\end{align*}
The linearization operator reads
\begin{align*}
&l_{\mathcal{E}} = \,\overline{\!D}_t - 12u_1^2\ \overline{\!D}_x - \,\overline{\!D}_x^{\,3}\,.
\end{align*}
Its adjoint operator has the form
\begin{align*}
&l_{\mathcal{E}}^{\,*} = -\,\overline{\!D}_t + 12u_1^2\ \overline{\!D}_x + 12\,\overline{\!D}_x(u_1^2) + \,\overline{\!D}_x^{\,3}\,.
\end{align*}
It cannot increase the order of a function by $3$ or more. Let us recall that cosymmetries of $\mathcal{E}$ are elements of $\ker l_{\mathcal{E}}^{\,*}$.

The equation $\mathcal{E}$ possesses the scaling symmetry
\begin{align*}
X = 3t\partial_t + x\partial_x - \sum_{j\,=\,0}^{+\infty} ju_j\partial_{u_j}\,.
\end{align*}
Note that $X$ cannot increase the order of a function. The characteristic $\varphi$ of $X$ has the form
\begin{align*}
\varphi = -3t(4u_1^3 + u_3) - xu_1\,.
\end{align*} 
In terms of cosymmetries, the Lie derivative $\mathcal{L}_X\colon E_1^{\hspace{0.1ex} 1,\hspace{0.2ex} 1}(\mathcal{E})\to E_1^{\hspace{0.1ex} 1,\hspace{0.2ex} 1}(\mathcal{E})$ is given by the operator $E_{\varphi} + l_{\varphi}^{\,*}$ or by $X + 1$. These operators coincide on cosymmetries.

The main goal of this section is to prove that the presymplectic structure $\Omega$ is not $d_1$-exact, i.e., that the element of $E_2^{\hspace{0.1ex} 2,\hspace{0.2ex} 1}(\mathcal{E})$ represented by $\Omega$ is non-trivial. In more applied terms, we need to show that there is no cosymmetry $\psi_0\in \ker l_{\mathcal{E}}^{\, *}$ such that $l_{\psi_0} - l_{\psi_0}^{\,*} = \Delta$.
For our further investigation, it is important that the presymplectic structure $\Omega$ is $X$-invariant, i.e., that $\mathcal{L}_X \Omega = 0$. This follows from the identity $l_{\psi} - l_{\psi}^{\, *} = 0$ for the cosymmetry $\psi = \Delta(\varphi)$. In addition, it is convenient to use the commutator
\begin{align}
[\partial_{u_{j}}, \,\overline{\!D}_t] = \partial_{u_{j}}\big(\,\overline{\!D}_x^{\,i}(4u_1^3) + u_{i+3}\big)\partial_{u_{i}}
\label{commutwithDt}
\end{align}
and the simple combinatorial

\vspace{1ex}

\propositiona{\label{Prop1} For any integers $i \geqslant 1$, $j\geqslant 0$,
\begin{align*}
\partial_{u_{j}} \,\overline{\!D}_x^{\,i} = \sum_{r\,=\,0}^{\min\{i\hspace{0.1ex},\,j\}}     \binom{i}{r} \,\overline{\!D}_x^{\,i - r} \partial_{u_{j-r}}\,.
\end{align*}
}

\vspace{-1.5ex}

Another useful fact is given by the following

\vspace{1ex}

\propositiona{\label{Prop2} If $\psi\in \mathcal{F}(\mathcal{E})$ and $s \geqslant 1$ is an integer such that $\ord \psi \leqslant s$, then
\begin{align*}
\partial_{u_{s+2}}\,l_{\mathcal{E}}^{\,*}(\psi) = 3\,\overline{\!D}_x (\partial_{u_s}\psi)\,.
\end{align*}
}

\vspace{-2ex}

\noindent
\textbf{Proof.} The order of the expression $12u_1^2\ \overline{\!D}_x(\psi) + 12\,\overline{\!D}_x(u_1^2)\psi$ is strictly less than $s + 2$. Then
\begin{align*}
\partial_{u_{s+2}}\,l_{\mathcal{E}}^{\,*}(\psi) = -\partial_{u_{s+2}} \,\overline{\!D}_t(\psi) + \partial_{u_{s+2}} \,\overline{\!D}_x^{\,3}(\psi)
\end{align*}
From Proposition~\ref{Prop1}, one obtains the formula $\partial_{u_{s+2}}\,\overline{\!D}_x^{\,3} = \,\overline{\!D}_x^{\,3} \partial_{u_{s+2}} + 3\,\overline{\!D}_x^{\,2} \partial_{u_{s+1}} + 3\,\overline{\!D}_x\hspace{0.15ex} \partial_{u_{s}} + \partial_{u_{s-1}}$.
Using~\eqref{commutwithDt}, we find
\begin{align*}
\partial_{u_{s+2}}\,l_{\mathcal{E}}^{\,*}(\psi) = - \partial_{u_{s+2}}\big(\,\overline{\!D}_x^{\,i}(4u_1^3) + u_{i+3}\big)\partial_{u_{i}}\psi + 3\,\overline{\!D}_x (\partial_{u_{s}}\psi) + \partial_{u_{s-1}}\psi = 3\,\overline{\!D}_x (\partial_{u_{s}}\psi)\,.
\end{align*}

\vspace{1ex}

In order to prove Theorem~\ref{MainTheo} by contradiction, we now assume that the presymplectic operator $\Delta = \,\overline{\!D}_x$ is produced by a cosymmetry. Let $\psi_0$ denote a cosymmetry such that $l_{\psi_0} - l_{\psi_0}^{\,*} = \Delta$, and let $k$ denote its order, $\ord \psi_0 = k$. \emph{We adopt this notation for all lemmas given in this section.} Direct computations of cosymmetries of orders $\leqslant 6$ show that 
$$
k \geqslant 7\,.
$$

\lemmaa{\label{Lemma1} There exist functions $a = a(t)$ and $\psi_2\in \mathcal{F}(\mathcal{E})$ such that
$\psi_0 = a u_k + \psi_2$, $\ord \psi_2\leqslant k-2$. Besides, $k$ is even.}

\vspace{0.5ex}

\noindent
\textbf{Proof.} Since $l_{\mathcal{E}}^{\,*}(\psi_0) = 0$, taking $s = k$ in Proposition~\ref{Prop2}, we get $\,\overline{\!D}_x (\partial_{u_{k}}\psi_0) = 0$. Then there exist functions $a = a(t)$ and $\psi_1\in \mathcal{F}(\mathcal{E})$ such that $\psi_0 = a u_{k} + \psi_1$, where $\ord \psi_1\leqslant k-1$. Note that $a$ cannot be identically zero (otherwise, $\ord \psi_0 < k = \ord \psi_0$, where $\ord \psi_0 \geqslant 7$).

The condition $l_{\psi_0} - l_{\psi_0}^{\,*} = \,\overline{\!D}_x$ implies that the operator 
$$
l_{\psi_0} - l_{\psi_0}^{\,*} = (1 - (-1)^k)\hspace{0.15ex} a \,\overline{\!D}_x^{\,k} + (1 - (-1)^{k-1})\hspace{0.15ex} (\partial_{u_{k-1}} \psi_1) \, \overline{\!D}_x^{\,k-1} + \ldots \,\overline{\!D}_x^{\,k-2} + \ldots
$$
can involve neither $\,\overline{\!D}_x^{\,k}$ nor $\,\overline{\!D}_x^{\,k-1}$. Then $k$ is even and $\partial_{u_{k-1}} \psi_1 = 0$. Consequently, $\psi_0$ has the form $\psi_0 = a u_{k} + \psi_2$, where $\ord \psi_2\leqslant k-2$.

\vspace{1ex}

Using Lemma~\ref{Lemma1}, one can derive a more precise form of $\psi_0$. In fact, Lemma~\ref{Lemma2}, Lemma~\ref{Lemma3}, and Lemma~\ref{Lemma4} are needed only to demonstrate that the coefficient of $u_k$, denoted throughout this section by $a$, is a constant.

\vspace{1ex}

\lemmaa{\label{Lemma2} There exist functions $b = b(t)$ and $\psi_3\in\mathcal{F}(\mathcal{E})$ such that $\psi_0 = a u_k + Bu_{k-2} + \psi_3$, where $\ord \psi_3\leqslant k-3$ and
\begin{align}
\label{BfromLemma2}
B = \dfrac{\dot{a}}{3}\hspace{0.1ex} x + 4(k-1)a u_1^2 + b\,.
\end{align}
}

\vspace{-2ex}

\noindent
\textbf{Proof.} According to Lemma~\ref{Lemma1}, $\psi_0 = a u_{k} + \psi_2$, where $\ord \psi_2\leqslant k-2$. Taking $s = k - 2$ in Proposition~\ref{Prop2}, we get the relation $\partial_{u_k} l_{\mathcal{E}}^{\,*}(\psi_2) = 3\,\overline{\!D}_x (\partial_{u_{k-2}}\psi_2)$. One can see that
\begin{align}
\label{firstterm}
-l_{\mathcal{E}}^{\,*}(au_k) = \dot{a}u_k + a\,\overline{\!D}_x^{\,k}(4u_1^3) - 12au_1^2u_{k+1} - 12a\,\overline{\!D}_x(u_1^2)u_k\,.
\end{align}
Using Proposition~\ref{Prop1}, one obtains $\partial_{u_k} \,\overline{\!D}_x^{\,k}(4u_1^3) = \binom{k}{k-1}\,\overline{\!D}_x\hspace{0.15ex} \partial_{u_1}(4u_1^3) = 12k \,\overline{\!D}_x(u_1^2)$. The condition $l_{\mathcal{E}}^{\,*}(\psi_0) = 0$ means that $l_{\mathcal{E}}^{\,*}(\psi_2) = -l_{\mathcal{E}}^{\,*}(au_k)$, and hence, $\partial_{u_k} l_{\mathcal{E}}^{\,*}(\psi_2) = -\partial_{u_k} l_{\mathcal{E}}^{\,*}(au_k)$. The relation
\begin{align*}
3\,\overline{\!D}_x (\partial_{u_{k-2}}\psi_2) = - \partial_{u_{k}} l_{\mathcal{E}}^{\,*}(au_k) = \dot{a} + 12(k-1)a \,\overline{\!D}_x(u_1^2)
\end{align*}
and the inequality $k - 2\geqslant 5$ complete the proof.

\vspace{1ex}

Now, in Lemma~\ref{Lemma3}, Lemma~\ref{Lemma4}, and Lemma~\ref{Lemma5}, we denote by $B$ the expression defined by~\eqref{BfromLemma2}.

\vspace{1ex}

\lemmaa{\label{Lemma3} There exists a function $\psi_4\in \mathcal{F}(\mathcal{E})$ such that $\psi_0 = a u_k + Bu_{k-2} + Cu_{k-3} + \psi_4$, where $\ord \psi_4\leqslant k-4$ and 
\begin{align*}
C = \dfrac{k-2}{2}\hspace{0.5ex} \overline{\!D}_x(B)\,.
\end{align*}
}

\vspace{-2ex}

\noindent
\textbf{Proof.} Since $k - 3\geqslant 4$ is odd, the condition $l_{\psi_0} - l_{\psi_0}^{\,*} = \,\overline{\!D}_x$ and Lemma~\ref{Lemma2} imply that the operator 
$$
l_{\psi_0} - l_{\psi_0}^{\,*} = \big(2\partial_{u_{k-3}}\psi_3 - (k-2)\,\overline{\!D}_x(B)\big)\,\overline{\!D}_x^{\,k-3} + \ldots \,\overline{\!D}_x^{\,k-4} + \ldots
$$
cannot involve $\,\overline{\!D}_x^{\,k-3}$. Then $2\partial_{u_{k-3}}\psi_3 = (k-2)\,\overline{\!D}_x(B)$. This observation completes the proof.

\vspace{1.5ex}

To make the proof of Lemma~\ref{Lemma5} clearer, we provide the following technical

\vspace{1.5ex}

\lemmaa{\label{Lemma4} There exists a function $f\in \mathcal{F}(\mathcal{E})$ such that $-\partial_{u_{k-2}}\,l_{\mathcal{E}}^{\,*}(Bu_{k-2}) = 4\dot{a}u_1^2 + \,\overline{\!D}_x(f)$.}

\vspace{0.5ex}

\noindent
\textbf{Proof.} Using Proposition~\ref{Prop1}, we get $\partial_{u_{k-2}}\,\overline{\!D}_x^{\,k-2}(4u_1^3 + u_3) = \binom{k-2}{k-3}\,\overline{\!D}_x\hspace{0.15ex} \partial_{u_1} (4u_1^3) = 12(k-2)\,\overline{\!D}_x (u_1^2)$. Since $k-2\geqslant 5 > 4 = \ord B + 3$, we obtain
\begin{align*}
-\partial_{u_{k-2}}\,l_{\mathcal{E}}^{\,*}(Bu_{k-2}) =&\ \,\overline{\!D}_t(B) + B\partial_{u_{k-2}}\,\overline{\!D}_x^{\,k-2}(4u_1^3 + u_3) - \,\overline{\!D}_x^{\,3}(B) - 12u_1^2\,\overline{\!D}_x(B) - 12\,\overline{\!D}_x(u_1^2)B =\\
&\ \,\overline{\!D}_t(B) + 12(k-2)B \,\overline{\!D}_x(u_1^2) - \,\overline{\!D}_x\big(\,\overline{\!D}_x^{\,2}(B) + 12u_1^2 B \big)\,.
\end{align*}
Then for $f_1 = -\,\overline{\!D}_x^{\,2}(B) + 12(k-3)u_1^2 B$, one has
\begin{align}
\notag
&-\partial_{u_{k-2}}\,l_{\mathcal{E}}^{\,*}(Bu_{k-2}) - \,\overline{\!D}_x(f_1) = \,\overline{\!D}_t(B) - 12(k-2)u_1^2 \,\overline{\!D}_x(B) ={}\\
&\dfrac{\ddot{a}}{3}\hspace{0.1ex} x + 12(k-1)\dfrac{\dot{a}}{3}\hspace{0.1ex} u_1^2 + 8(k-1)a u_1\,\overline{\!D}_x(4u_1^3 + u_3) + \dot{b} - 12(k-2)u_1^2 \Big( \dfrac{\dot{a}}{3}\hspace{0.1ex} + 4(k-1)a \,\overline{\!D}_x(u_1^2) \Big).
\label{xtotalder}
\end{align}
Each summand in~\eqref{xtotalder} belongs to the image of $\,\overline{\!D}_x$, except for the terms that contain $\dot{a}$. The sum of these terms is $4\dot{a}u_1^2$. This observation completes the proof.

\vspace{1ex}

Finally, we arrive at the following result

\vspace{1ex}

\lemmaa{\label{Lemma5} There exist $\psi_1\in \mathcal{F}(\mathcal{E})$ and $a_0\in \mathbb{R}$ such that $\psi_0 = a_0 u_k + \psi_1$, $a_0\neq 0$, and $\ord \psi_1\leqslant k-1$.}

\vspace{0.5ex}

\noindent
\textbf{Proof.} The condition $\partial_{u_{k-2}}l_{\mathcal{E}}^{\,*}(\psi_0) = 0$ and Lemma~\ref{Lemma3} give rise to the relation 
$$
\partial_{u_{k-2}}l_{\mathcal{E}}^{\,*}(\psi_4) = -\partial_{u_{k-2}}l_{\mathcal{E}}^{\,*}(au_k + Bu_{k-2} + Cu_{k-3})\,.
$$
Since $\ord \psi_4\leqslant k-4\geqslant 3$, from Proposition~\ref{Prop2} it follows that the LHS belongs to the image of $\,\overline{\!D}_x$. Therefore, the RHS also belongs to the image of $\,\overline{\!D}_x$. Using~\eqref{firstterm} and Proposition~\ref{Prop1}, one can see that $-\partial_{u_{k-2}}l_{\mathcal{E}}^{\,*}(au_k)\in \mathrm{im}\, \,\overline{\!D}_x$. Next, by Lemma~\ref{Lemma1}, $k-2\geqslant 6$, and a simple calculation shows that
\begin{align*}
-\partial_{u_{k-2}}l_{\mathcal{E}}^{\,*}(Cu_{k-3}) = -3\,\overline{\!D}_x^{\, 2}(C)\in \mathrm{im}\, \,\overline{\!D}_x\,.
\end{align*}
Here $\ord C = 2$.
Ultimately, from Lemma~\ref{Lemma4}, one sees that $\dot{a} = 0$ and hence $a = a_0$ is a constant.

\vspace{1ex}

Of course, we now know more about these $\psi_1$ and $k$, but at this point, the additional information turns out to be redundant. The final step is given by

\vspace{1ex}

\theorema{\label{MainTheo} The element of $E_2^{\hspace{0.1ex} 2,\hspace{0.2ex} 1}(\mathcal{E})$ represented by the presymplectic structure $\Omega$ is non-trivial.}

\vspace{0.5ex}

\noindent
\textbf{Proof.} Let $\psi$ be a lowest order cosymmetry that produces $\Delta$. Denote $\ord \psi = k$. Then $k \geqslant 7$ and, according to Lemma~\ref{Lemma5}, $\psi = a_0 u_k + \psi_1$, where $\ord \psi_1 \leqslant k-1$ and $a_0\in\mathbb{R}\setminus\{0\}$. Let $\omega_{\psi}$ denote the variational $1$-form that corresponds to $\psi$. Then the cosymmetry of the Lie derivative $\mathcal{L}_X \omega_{\psi}$ is
\begin{align}
(X+1)\psi = (-k+1)a_0 u_k + (X+1)\psi_1 = (-k+1)\psi + (X+k)\psi_1\,.
\label{Liedercosym}
\end{align}
Note that $(X+k)\psi_1 = (X+k)\psi$ is a cosymmetry, and its order is strictly less than $k$. Denote by $\widetilde{\omega}$ the variational $1$-form that corresponds to $(X+k)\psi_1$. The formula~\eqref{Liedercosym} results in
\begin{align*}
\mathcal{L}_X \omega_{\psi} = (-k+1) \omega_{\psi} + \widetilde{\omega}\,.
\end{align*}
Since $\Omega = d_1 \omega_{\psi}$ and $\mathcal{L}_X \Omega = 0$, applying $d_1$, we obtain
\begin{align*}
0 = (-k+1) \Omega + d_1 \widetilde{\omega}\,.
\end{align*}
Here $-k+1 \neq 0$. Thus the variational $1$-form
\begin{align*}
\dfrac{1}{k-1}\, \widetilde{\omega}
\end{align*}
is also a potential for $\Omega$. However, the order of its cosymmetry is strictly less than $k$.
This contradiction completes the proof.

\section{\label{Section4}On conservation laws with trivial cosymmetries}

Any $\ell$-normal system of differential equations has trivial group $\ker d_1^{\hspace{0.1ex} 0, \hspace{0.2ex} n-1} = E_2^{\hspace{0.1ex} 0, \hspace{0.2ex} n-1}$ of the Vinogradov $\mathcal{C}$-spectral sequence~\cite{VinKr}. This means that any non-trivial conservation law of an $\ell$-normal system corresponds to a non-trivial cosymmetry. The same conclusion holds for any system that admits a symmetry group scaling each dependent variable by a factor with a positive exponent while preserving all independent variables (see Proposition on p.~111 in~\cite{Vin} and Theorem~12 in~\cite{Olver0} for the case of equal exponents). We propose Example~\ref{MainExamCL}, where a non-trivial conservation law corresponds to a trivial cosymmetry. The characteristic described there vanishes on the corresponding (overdetermined) system.

A significant observation is given by

\vspace{1ex}

\theorema{\label{Theorem2} Let $L\in \Lambda_h^n(\pi)$ be a horizontal $n$-form such that the variational derivative $\mathrm{E}(L)$ vanishes on an infinitely prolonged system $\mathcal{E}\subset J^{\infty}(\pi)$. If the element of $E_{2}^{\hspace{0.1ex} 2, \hspace{0.2ex} n-1}(\mathcal{E})$ represented by the corresponding presymplectic structure of $\mathcal{E}$ is non-trivial, then $L|_{\mathcal{E}}\in E_{0}^{\hspace{0.1ex} 0, \hspace{0.2ex} n}(\mathcal{E})$ represents a non-trivial element of $E_{1}^{\hspace{0.1ex} 0, \hspace{0.2ex} n}(\mathcal{E})$, which represents a non-trivial element of $E_{2}^{\hspace{0.1ex} 0, \hspace{0.2ex} n}(\mathcal{E})$.}

\vspace{0.5ex}

\noindent
\textbf{Proof.} Let $l\in \Lambda^n(\mathcal{E})$ be a differential form that represents the internal Lagrangian determined by $L$. Then $dl$ gives rise to the presymplectic structure, and $l \in L|_{\mathcal{E}}$. Since the differential
\begin{align*}
d_2^{\hspace{0.2ex} 0,\hspace{0.2ex} n}\colon E_2^{\hspace{0.1ex} 0, \hspace{0.2ex} n}(\mathcal{E}) \to E_2^{\hspace{0.1ex} 2, \hspace{0.2ex} n-1}(\mathcal{E})
\end{align*} 
is induced by the de Rham differential $d$, the non-triviality of the corresponding element of $E_{2}^{\hspace{0.1ex} 2, \hspace{0.2ex} n-1}(\mathcal{E})$ implies that the element of $E_2^{\hspace{0.1ex} 0, \hspace{0.2ex} n}(\mathcal{E})$ determined by $L|_{\mathcal{E}}$ cannot be zero.

\vspace{1ex}

This theorem shows, for example, that the restriction of $\lambda$ from~\eqref{almostpKdVLagr} to the infinite prolongation of~\eqref{almostpKdV} is not a total divergence on the potential mKdV. 

\vspace{1ex}

\remarka{Even though some presymplectic structures coming from Lagrangians are not produced by variational $1$-forms, they are produced by internal Lagrangians. Let us note that for such a presymplectic structure and any spatial distribution $\mathcal{S}$ (see~\cite{Druzhkov4}, Section~3) of a differential equation, the corresponding $\mathcal{S}$-presymplectic structure still comes from an $\mathcal{S}$-variational $1$-form (Remark~3 in~\cite{Druzhkov4}). There is some difference between the covariant and non-covariant approaches in this case.}

\vspace{1ex}

The following trick allows one to formally increase the number of independent variables by~$1$ and reinterpret~\eqref{almostpKdVLagr} in terms of conservation laws. The vector field $\partial_y$ on $\mathbb{R}$ (with a global coordinate $y$) makes it a differential equation. We denote it by $\mathbb{R}_{\partial_y}$. Let $\mathcal{E}$ be a differential equation. Denote by $\mathcal{E}'_{\partial_y}$ the product 
$$
\mathcal{E}'_{\partial_y} = \mathcal{E}\times \mathbb{R}_{\partial_y}\,.
$$
The plane of the Cartan distribution of $\mathcal{E}'_{\partial_y}$ at a point $(\rho, y)$ is spanned by the plane of the Cartan distribution of $\mathcal{E}$ at $\rho\in \mathcal{E}$ and the vector of $\partial_y$ at $y\in \mathbb{R}_{\partial_y}$. In simpler terms, compared to $\mathcal{E}$, the equation $\mathcal{E}'_{\partial_y}$ has the extra independent variable $y$ and the additional total derivative $\partial_y$. 

\vspace{1ex}

\examplea{\label{MainExamCL} Let $\mathcal{E}$ again denote the infinite prolongation of the potential mKdV~\eqref{almostpKdV}. 
The equation
$\mathcal{E}'_{\partial_y}$ can be described as the infinite prolongation of the overdetermined system
\begin{align}
&u_t = 4u_x^3 + u_{xxx}\,,\qquad u_y = 0
\label{almostpKdV_yCl}
\end{align}
with the independent variables $t$, $x$, $y$. Consider its conservation law given by the total divergence
\begin{align}
D_t(0) + D_x(0) + D_y(\lambda) = 0\,, \qquad \lambda = \dfrac{u_t u_x}{2} - u_x^4 + \dfrac{u_{xx}^2}{2}\,.
\label{TheCL}
\end{align}
According to the decomposition
\begin{align*}
D_y(\lambda) = u_{xy}\hspace{0.2ex} (u_t - 4u_x^3 - u_{xxx}) + 0\!\hspace{0.15ex}\cdot\!\hspace{0.15ex} u_y + D_t\Big(\dfrac{u_x}{2}\hspace{0.15ex} u_y\Big) + D_x\Big(u_{xx} u_{xy} - \dfrac{u_t}{2}\hspace{0.15ex} u_y\Big),
\end{align*}
the conservation law is associated, for example, with the characteristic $Q$
given by
\begin{align*}
Q_1 = u_{xy}\,,\qquad Q_2 = 0\,.
\end{align*}
Then the cosymmetry $Q|_{\mathcal{E}'_{\partial_y}}$ is zero. Assume that the conservation law is trivial. Then there exist functions $g_1, g_2 \in \mathcal{F}(\mathcal{E}'_{\partial_y})$ such that
\begin{align*}
\lambda|_{\mathcal{E}'_{\partial_y}} = \widetilde{D}_t(g_2) - \widetilde{D}_x(g_1)\,, \qquad \widetilde{D}_x = D_x|_{\mathcal{E}'_{\partial_y}},\ \
\widetilde{D}_t = D_t|_{\mathcal{E}'_{\partial_y}}\,.
\end{align*}
Abusing the notation, one sees that in this case, $\lambda|_{\mathcal{E}}$ is a total divergence on the potential mKdV, $\lambda|_{\mathcal{E}} = \,\overline{\!D}_t(g_2|_{y\,=\,0}) - \,\overline{\!D}_x(g_1|_{y\,=\,0})$. This relation contradicts to Theorem~\ref{Theorem2} applied to the Lagrangian~\eqref{almostpKdVLagr}. Consequently, the conservation law~\eqref{TheCL} of~\eqref{almostpKdV_yCl} is non-trivial.

\vspace{1ex}

\remarka{The horizontal form corresponding to
$$
u_x^4\hspace{0.2ex} dt\wedge dx\in \Lambda^2(\mathcal{E}'_{\partial_y})
$$
represents the conservation law~\eqref{TheCL} since it coincides with the restriction of 
\begin{align*}
\lambda\, dt\wedge dx +  d_h\Big(\dfrac{u_x u_{xx}}{2}\, dt\Big)
\end{align*}
to the system $\mathcal{E}'_{\partial_y}$.
}
}

\subsection{Homotopy of differential equations}

A more geometric approach is based on the notion of homotopy of differential equations~\cite{Vin} (see p.~111). Using it, one can show that the projection onto the first factor $\mathcal{E}\times \mathbb{R}_{\partial_y} \to \mathcal{E}$
induces an isomorphism of $E_2(\mathcal{E})$ and $E_2(\mathcal{E}'_{\partial_y})$. Recall the details.
Let $\mathcal{E}_i$ be differential equations with Cartan distributions $\mathcal{C}_i$, $i = 1, 2$. Denote by $\mathcal{C}_{i\, \rho}$ the Cartan plane at a point $\rho\in \mathcal{E}_i$.

\vspace{1ex}

\definitiona{A mapping $f\colon \mathcal{E}_1\to \mathcal{E}_2$ is a DE-morphism if its differential maps Cartan planes of $\mathcal{E}_1$ to the respective Cartan planes of $\mathcal{E}_2$, i.e., for each point $\rho\in \mathcal{E}_1$, 
\begin{align*}
df_\rho(\mathcal{C}_{1\, \rho})\subset \mathcal{C}_{2\, f(\rho)}.
\end{align*}
}

\noindent
In particular, point and contact transformations of differential equations are DE-morphisms, as well as any differential coverings~\cite{VinKr}.

If $f\colon \mathcal{E}_1\to \mathcal{E}_2$ is a DE-morphism, its pullback $f^*\colon \Lambda^*(\mathcal{E}_2) \to \Lambda^*(\mathcal{E}_1)$ induces a homomorphism of the Cartan ideals, i.e., $f^*(\mathcal{C}\Lambda^*(\mathcal{E}_2)) \subset \mathcal{C}\Lambda^*(\mathcal{E}_1)$; accordingly, $f^*(\mathcal{C}^p\Lambda^{p+q}(\mathcal{E}_2)) \subset \mathcal{C}^p\Lambda^{p+q}(\mathcal{E}_1)$. Therefore, DE-morphisms of differential equations induce homomorphisms of their $\mathcal{C}$-spectral sequences.

\vspace{1ex}

\definitiona{DE-morphisms $f_0, f_1\colon \mathcal{E}_1\to \mathcal{E}_2$ are homotopic if there is a mapping $f\colon \mathcal{E}_1\times [0; 1]\to \mathcal{E}_2$ such that $f(\rho, 0) = f_0(\rho)$ and $f(\rho, 1) = f_1(\rho)$ for each $\rho\in\mathcal{E}_1$, and all the mappings $f_\tau\colon \mathcal{E}_1\to \mathcal{E}_2$, $f_{\tau}\colon \rho\mapsto f(\rho, \tau)$ are DE-morphisms (for $\tau\in [0; 1]$).
}

\vspace{1ex}

Homotopic DE-morphisms of equations induce the same homomorphisms of the second pages of their $\mathcal{C}$-spectral sequences~\cite{Vin} (see also Appendix~\ref{App:B}).
The projection onto the first factor
\begin{align*}
\mathrm{pr}_{\mathcal{E}} \colon \mathcal{E}\times \mathbb{R}_{\partial_y}\to \mathcal{E}
\end{align*}
and the mapping $\mathcal{E}\to \mathcal{E}\times \mathbb{R}_{\partial_y}$, $\rho\mapsto (\rho, 0)$ establish a homotopy equivalence of $\mathcal{E}$ and $\mathcal{E}'_{\partial_y}$. Then we obtain the following

\vspace{1ex}

\theorema{Let $L\in \Lambda_h^n(\pi)$ be a horizontal $n$-form such that the variational derivative $\mathrm{E}(L)$ vanishes on an infinitely prolonged system $\mathcal{E}\subset J^{\infty}(\pi)$. If the element of $E_{2}^{\hspace{0.1ex} 2, \hspace{0.2ex} n-1}(\mathcal{E})$ represented by the corresponding presymplectic structure of $\mathcal{E}$ is non-trivial, then $\mathrm{pr}_{\mathcal{E}}^{\,*}(L|_{\mathcal{E}})\in E_{0}^{\hspace{0.1ex} 0, \hspace{0.2ex} n} (\mathcal{E}'_{\partial_y})$ represents a non-trivial conservation law of $\mathcal{E}'_{\partial_y}$ with a trivial cosymmetry.
}

\vspace{1ex}

Let us recall that the Noether correspondence between symmetries and variational $1$-forms of a Lagrangian equation is surjective (roughly, in this case, cosymmetries are symmetries). We get

\vspace{1ex}

\theorema{\label{TheorLagrExa}Let $\mathcal{E}$ be the infinite prolongation of a system of Euler-Lagrange equations $\mathrm{E}(L) = 0$ that possesses only Noether symmetries (or no symmetries). If the corresponding presymplectic structure of $\mathcal{E}$ is non-trivial, then the system $\mathcal{E}'_{\partial_y}$ admits a non-trivial conservation law with a trivial cosymmetry.
}

\vspace{1ex}

\remarka{Any differential equation is homotopic to its tangent equation~\cite{KraVer1}.
}

\subsection{Example~\ref{MainExamCL} as a symmetry reduction}

System~\eqref{almostpKdV_yCl} can be interpreted as a symmetry reduction. For instance, the infinite prolongation of the equation
\begin{align}
u_t = 4u_x^3 + u_{xxx} + u_{yyy}
\label{reduction}
\end{align}
has the symmetry $\partial_y$. The corresponding invariant solutions are described by~\eqref{almostpKdV_yCl}.

\vspace{1ex}

\remarka{If the variational derivative $\mathrm{E}(L)$ of a Lagrangian $L\in \Lambda_h^n(\pi)$ vanishes on an infinitely prolonged system $\mathcal{N}$, the Noether theorem relates some symmetries of $\mathcal{N}$ and its conservation laws. The relation is based on the formula~\eqref{Noethiden} and can be described in terms of the corresponding presymplectic structure $\Omega_L$. Namely, a symmetry $Y$ of $\mathcal{N}$ is a Noether symmetry if the variational $1$-form $Y\lrcorner\, \Omega_L$ is $d_1$-exact. Then $Y$ corresponds to a conservation law $\xi\in E_1^{\hspace{0.1ex} 0, \hspace{0.2ex} n-1}(\mathcal{N})$ such that
\begin{align}
Y\lrcorner\, \Omega_L = d_1\xi\,.
\label{Noethercorr}
\end{align}
Denote by $\mathcal{N}_Y$ the system that describes $Y$-invariant solutions of $\mathcal{N}$, $\mathcal{N}_Y\subset \mathcal{N}$. The restriction of $Y$ to $\mathcal{N}_Y$ is a trivial symmetry. Then, restricting~\eqref{Noethercorr} to $\mathcal{N}_Y$, we get $d_1(\xi|_{\mathcal{N}_Y}) = 0$. Thus $\xi|_{\mathcal{N}_Y}$ is an element of $E_2^{\hspace{0.1ex} 0, \hspace{0.2ex} n-1}(\mathcal{N}_Y)$. The example~\eqref{reduction} shows that the conservation law $\xi|_{\mathcal{N}_Y}$ of $\mathcal{N}_Y$ is not necessarily trivial. Indeed, the Lagrangian
\begin{align*}
\Big( \dfrac{u_t u_x}{2} - u_x^4 + \dfrac{u_{xx}^2}{2} + \dfrac{u_{xy}u_{yy}}{2} \Big)dt\wedge dx\wedge dy
\end{align*}
produces a presymplectic structure for~\eqref{reduction}. Let us set it as $\Omega_L$. The symmetry $Y = -\partial_y$ is a Noether symmetry,
\begin{align}
-\partial_y \hspace{0.1ex} \lrcorner\, \Omega_L = d_1 \xi\,,
\label{Noetherred}
\end{align}
where the conservation law $\xi$ is defined by the restriction of
\begin{align*}
- \dfrac{u_x u_y}{2}\hspace{0.15ex} dx\wedge dy + \Big(u_{xx} u_{xy} - \dfrac{u_t}{2}\hspace{0.15ex} u_y - \dfrac{u_{yy}^2}{2}\Big) dt\wedge dy + (\lambda - u_{xy} u_{yy})dt\wedge dx\,,\quad \lambda = \dfrac{u_t u_x}{2} - u_x^4 + \dfrac{u_{xx}^2}{2}
\end{align*}
to the infinite prolongation of~\eqref{reduction}, which we take as $\mathcal{N}$. The corresponding system $\mathcal{N}_Y$ coincides with the infinite prolongation of~\eqref{almostpKdV_yCl}.
Restricting $\xi|_{\mathcal{N}}$ to $\mathcal{N}_Y$, we get the non-trivial conservation law defined by~\eqref{TheCL}.

}

\section{Conclusion}

The results of this paper confirm that, in general, the relationship between proper conservation laws and their cosymmetries is many-to-many. An interesting but difficult direction for future research, in fact suggested in a different form by P.~Olver in~\cite{Olver0}, is to find a Lagrangian system whose group $E_2^{\hspace{0.1ex} 0, \hspace{0.2ex} n-1}$ has non-trivial, non-topological elements. In this case, Noether's theorem would fail to describe all proper conservation laws. Such a system would necessarily be a gauge system, meaning that there exists a non-trivial total differential operator whose image consists of characteristics of symmetries restricted to the system. Moreover, the system must not admit any symmetry group scaling each dependent variable by a factor with a positive exponent while preserving all independent variables.

Another intriguing problem is to construct an example of a non-topological, hidden internal Lagrangian, i.e., one whose corresponding presymplectic structure is zero. This problem is analogous to the one solved in this paper but seems to be significantly more challenging. 

\vspace{-1ex}

\section*{Acknowledgments}

The author appreciates valuable discussions with A. Shevyakov, J. Krasil'shchik, A. Verbovetsky, and P. Olver.
The author is also grateful to the University of Saskatchewan for its hospitality. Special thanks go to Prof. Alexey Shevyakov for financial support through the NSERC grant RGPIN 04308-2024 and to the Pacific Institute for the Mathematical Sciences for support through a PIMS Postdoctoral Fellowship.


\appendix

\section{Cosymmetries and variational $1$-forms}\label{App:A}

It is instructive to recall why the map from cosymmetries to variational $1$-forms is surjective and has a non-trivial kernel if the system is not $\ell$-normal. A more conceptual approach is based on the compatibility complex for $l_{\mathcal{E}}$~\cite{KraVer2} (see Corollary 7.4).

Denote by $\mathcal{C}(\varkappa(\pi), \Lambda^k_h(\pi))$ the $\mathcal{F}(\pi)$-module of $\mathcal{C}$-differential operators from $\varkappa(\pi)$ to $\Lambda^k_h(\pi)$. Let $\mathcal{C}(\varkappa(\mathcal{E}), \Lambda^{k}_h(\mathcal{E}))$ be the $\mathcal{F}(\mathcal{E})$-module of the restrictions of such operators to a system $\mathcal{E}$. Then
\begin{align*}
\mathcal{C}(\varkappa(\mathcal{E}), \Lambda^{k}_h(\mathcal{E})) = \mathcal{C}(\varkappa(\pi), \Lambda^k_h(\pi))/ I\cdot \mathcal{C}(\varkappa(\pi), \Lambda^k_h(\pi))\,,
\end{align*}
where $I\subset \mathcal{F}(\pi)$ is the ideal of $\mathcal{E}$. Analogously, we introduce $\mathcal{C}(P(\mathcal{E}), \Lambda^{k}_h(\mathcal{E}))$.
There is the canonical isomorphism 
\begin{align*}
\mathcal{C}(\varkappa(\mathcal{E}), \mathcal{F}(\mathcal{E})) = \mathcal{C}\Lambda^1(\pi)/I \cdot \mathcal{C}\Lambda^1(\pi)\,,
\end{align*}
which maps an operator $A|_{\mathcal{E}}$ to the coset $\omega_{A} + I \cdot \mathcal{C}\Lambda^1(\pi)$, where $\omega_{A}\in \mathcal{C}\Lambda^1(\pi)$ is a unique Cartan $1$-form such that for $\varphi\in \varkappa(\pi)$, $A(\varphi) = E_{\varphi}\lrcorner\, \omega_{A}$. In adapted coordinates, for $A(\varphi) = A^{\alpha}_iD_{\alpha}(\varphi^i)$, one has $\omega_A = A^{\alpha}_i\theta^i_{\alpha}$.
This isomorphism gives rise to the epimorphism 
$$
\nu^{0}\colon \mathcal{C}(\varkappa(\mathcal{E}), \mathcal{F}(\mathcal{E}))\to \mathcal{C}\Lambda^1(\mathcal{E})\,,\ \ A|_{\mathcal{E}}\mapsto \omega_A|_{\mathcal{E}}
$$
making the following sequence exact~\cite{VinKr} (see p.~198).
\begin{align*}
\xymatrix{
\mathcal{C}(P(\mathcal{E}), \mathcal{F}(\mathcal{E})) \ar[r]^-{\circ\hspace{0.15ex} l_{\mathcal{E}}} & \mathcal{C}(\varkappa(\mathcal{E}), \mathcal{F}(\mathcal{E})) \ar[r]^-{\nu^{0}} & \mathcal{C}\Lambda^1(\mathcal{E}) \ar[r] & 0
}
\end{align*}
This allows one to describe the horizontal de Rham complex $\Lambda_h^{\bullet}(\mathcal{E})$ with the desired coefficients.
Since $E_0^{\hspace{0.1ex} 1, \hspace{0.2ex} k}(\mathcal{E}) = \Lambda^{k}_h(\mathcal{E}) \otimes_{\mathcal{F}(\mathcal{E})} \mathcal{C}\Lambda^1(\mathcal{E})$ and $\mathcal{C}(\varkappa(\mathcal{E}), \Lambda^{k}_h(\mathcal{E})) = \Lambda^{k}_h(\mathcal{E}) \otimes_{\mathcal{F}(\mathcal{E})} \mathcal{C}(\varkappa(\mathcal{E}), \mathcal{F}(\mathcal{E}))$, we can set 
\begin{align*}
\nu^k\colon \mathcal{C}(\varkappa(\mathcal{E}), \Lambda^{k}_h(\mathcal{E})) \to E_0^{\hspace{0.1ex} 1, \hspace{0.2ex} k}(\mathcal{E})\quad \text{such that} \quad \nu^k\colon [\omega]_{h} \otimes A|_{\mathcal{E}} \mapsto \omega \wedge \nu^0(A|_{\mathcal{E}}) + \mathcal{C}^2\Lambda^{k+1}(\mathcal{E}) \,,
\end{align*}
where $\omega\in \Lambda^k(\mathcal{E})$, $[\omega]_h = \omega + \mathcal{C}\Lambda^k(\mathcal{E})$, and each $\nu^k$ is a homomorphism of the $\mathcal{F}(\mathcal{E})$-modules. For instance, if $B(\varphi) = B_{ij}^{\alpha}D_{\alpha}(\varphi^i)dx^j$ for $\varphi\in \varkappa(\pi)$, then $\nu^1(B|_{\mathcal{E}})\in E_0^{\hspace{0.1ex} 1, \hspace{0.2ex} 1}(\mathcal{E})$ is determined by the differential $2$-form $(B_{ij}^{\alpha}\hspace{0.15ex} dx^j\wedge \theta^i_{\alpha})|_{\mathcal{E}}$.
Then all columns of the commutative diagram
\begin{align*}
\xymatrix{
\ldots \ar[r] & \mathcal{C}(P(\mathcal{E}), \Lambda_h^{n-2}(\mathcal{E}))\ar[r]^-{d_h\circ } \ar[d]^{\circ\hspace{0.15ex} l_{\mathcal{E}}} & \mathcal{C}(P(\mathcal{E}), \Lambda_h^{n-1}(\mathcal{E}))\ar[r]^-{d_h\circ } \ar[d]^{\circ\hspace{0.15ex} l_{\mathcal{E}}} & \mathcal{C}(P(\mathcal{E}), \Lambda_h^{n}(\mathcal{E})) \ar[r] \ar[d]^{\circ\hspace{0.15ex} l_{\mathcal{E}}} & 0\\
\ldots \ar[r] & \mathcal{C}(\varkappa(\mathcal{E}), \Lambda_h^{n-2}(\mathcal{E}))\ar[r]^-{d_h\circ } \ar[d]^-{\nu^{n-2}} & \mathcal{C}(\varkappa(\mathcal{E}), \Lambda_h^{n-1}(\mathcal{E}))\ar[r]^-{d_h\circ } \ar[d]^-{\nu^{n-1}} & \mathcal{C}(\varkappa(\mathcal{E}), \Lambda_h^{n}(\mathcal{E})) \ar[r] \ar[d]^-{\nu^{n}} & 0\\
\ldots \ar[r] & E_0^{\hspace{0.1ex} 1, \hspace{0.2ex} n-2}(\mathcal{E}) \ar[r]^-{d_0} \ar[d] & E_0^{\hspace{0.1ex} 1, \hspace{0.2ex} n-1}(\mathcal{E}) \ar[r]^-{d_0} \ar[d] & E_0^{\hspace{0.1ex} 1, \hspace{0.2ex} n}(\mathcal{E}) \ar[r] \ar[d] & 0\\
& 0 & 0 & 0
}
\end{align*}
are exact. The only non-trivial cohomology of the second row from the top appears at the term $\mathcal{C}(\varkappa(\mathcal{E}), \Lambda_h^{n}(\mathcal{E}))$. It is given by $\widehat{\varkappa}(\mathcal{E})$~\cite{VinKr}.
Similarly, for the top row, the only non-trivial cohomology is given by $\widehat{P}(\mathcal{E})$. Recall that $\ker l_{\mathcal{E}}^{\, *} \subset \widehat{P}(\mathcal{E}) \subset \mathcal{C}(P(\mathcal{E}), \Lambda_h^{n}(\mathcal{E}))$. Diagram chasing shows that the map from $\ker l_{\mathcal{E}}^{\, *}$ to variational $1$-forms generated by $\psi\mapsto (-1)^{n-1} \nu^{n-1} ((d_h\hspace{0.1ex} \circ)^{\hspace{0.2ex} -1} (\psi\circ l_{\mathcal{E}}))$ is well-defined and surjective. This map coincides with the one determined by~\eqref{Cosymtoonef}. 

If $\psi = \nabla^*(\eta)$ for some $\eta\in\Lambda^n_h(\mathcal{E})$ and some $\mathcal{C}$-differential operator $\nabla\colon P(\mathcal{E})\to \mathcal{F}(\mathcal{E})$ such that $\nabla \circ l_{\mathcal{E}} = 0$, then $\psi$ is a cosymmetry, as $l_{\mathcal{E}}^{\,*}\circ \nabla^* = (\nabla \circ l_{\mathcal{E}})^* = 0$. According to Green's formula, $\psi\circ l_{\mathcal{E}}(\varphi) = \langle \nabla^*(\eta), l_{\mathcal{E}}(\varphi) \rangle = \langle \eta, \nabla\circ l_{\mathcal{E}}(\varphi) \rangle + d_h (\gamma\circ l_{\mathcal{E}}(\varphi)) = d_h (\gamma\circ l_{\mathcal{E}}(\varphi))$ for some operator $\gamma\in \mathcal{C}(P(\mathcal{E}), \Lambda_h^{n-1}(\mathcal{E}))$. Since $\gamma\circ l_{\mathcal{E}}\in \ker \nu^{n-1}$, this $\psi$ corresponds to the zero variational $1$-form.

\section{Homotopy and second pages of the $\mathcal{C}$-spectral sequence}\label{App:B}

\lemmaa{Let $f_0, f_1\colon \mathcal{E}_1\to \mathcal{E}_2$ be DE-morphisms of differential equations. If they are related by a homotopy $f\colon \mathcal{E}_1\times [0;1]\to \mathcal{E}_2$, then $f_0^*$ and $f_1^*$ induce the same homomorphisms $E^{\hspace{0.1ex} p,\hspace{0.15ex} q}_2(\mathcal{E}_2)\to E^{\hspace{0.1ex} p,\hspace{0.15ex} q}_2(\mathcal{E}_1)$.
}

\vspace{0.5ex}
\noindent
\textbf{Proof.} For each $\tau\in [0; 1]$, denote by $s_{\tau}$ the section $\mathcal{E}_1 \to \mathcal{E}_1\times [0;1]$, $\rho\mapsto (\rho, \tau)$. Then the homotopy formula
\begin{align*}
f_1^*(\omega) - f_0^*(\omega) = dK(\omega) + K(d\omega),\qquad K(\omega) = \int_0^1 s_\tau^*({\partial_{\tau}}\, \lrcorner\, f^*(\omega)) d\tau
\end{align*}
holds for any $\omega\in \Lambda^*(\mathcal{E}_2)$. At each point $\rho\in \mathcal{E}_1$, the relation
\begin{align*}
s_\tau^*({\partial_{\tau}}\hspace{0.2ex} \lrcorner\, f^*(\omega)) = f_{\tau}^*({df_{(\rho, \tau)}(\partial_{\tau})} \hspace{0.1ex} \lrcorner \,\omega)
\end{align*}
holds for any $\tau \in [0; 1]$. Since all $f_{\tau}$ are DE-morphisms, $K(\mathcal{C}^p\Lambda^{p+q}(\mathcal{E}_2))\subset \mathcal{C}^{p-1}\Lambda^{p-1+q}(\mathcal{E}_1)$ for $p>0$.
Now suppose that $\omega\in \mathcal{C}^p\Lambda^{p+q}(\mathcal{E}_2)$ gives rise to a $d_1$-closed element of $E^{\hspace{0.1ex} p,\hspace{0.15ex} q}_1(\mathcal{E}_2)$. One can describe the homomorphism of the first pages induced by $f^*_1 - f^*_0\colon \Lambda^*(\mathcal{E}_2) \to \Lambda^*(\mathcal{E}_1)$
using any representatives. Then we can assume, without loss of generality, that $d\omega \in \mathcal{C}^{p+2}\Lambda^{p+q+1}(\mathcal{E}_2)$. Due to the inclusion $K(\mathcal{C}^{p+2}\Lambda^{p+q+1}(\mathcal{E}_2))\subset \mathcal{C}^{p+1}\Lambda^{p+q}(\mathcal{E}_1)$, the differential form $K(d\omega)$ gives rise to the trivial element of $E_1^{\hspace{0.1ex} p,\hspace{0.2ex} q}(\mathcal{E}_1)$. Since $K(\omega)\in \mathcal{C}^{p-1}\Lambda^{p-1+q}(\mathcal{E}_1)$ and $dK(\omega)\in \mathcal{C}^{p}\Lambda^{p+q}(\mathcal{E}_1)$, the differential form $K(\omega)$ determines an element of $E_1^{\hspace{0.1ex} p-1,\hspace{0.2ex} q}(\mathcal{E}_1)$, and hence, $dK(\omega)$ leads to an element of $d_1 (E_1^{\hspace{0.1ex} p-1,\hspace{0.2ex} q}(\mathcal{E}_1))$. Thus the differential form $f_1^*(\omega) - f_0^*(\omega)$ determines an element of $d_1 (E_1^{\hspace{0.1ex} p-1,\hspace{0.2ex} q}(\mathcal{E}_1))$ and the homomorphism $E^{\hspace{0.1ex} p,\hspace{0.15ex} q}_2(\mathcal{E}_2)\to E^{\hspace{0.1ex} p,\hspace{0.15ex} q}_2(\mathcal{E}_1)$ induced by the mapping $f^*_1 - f^*_0$ is trivial. The case $p = 0$ differs in that the form $dK(\omega)$ results in the trivial element of $E_1^{\hspace{0.1ex}0, \hspace{0.2ex} q}(\mathcal{E}_1)$.

\end{document}